\documentclass[11pt, reqno]{article}

\usepackage[mathscr]{eucal}
\usepackage[all]{xy}
\usepackage{mathrsfs}
\usepackage{xypic}
\usepackage{amsfonts}
\usepackage{amsmath}
\usepackage{amsthm}
\usepackage{amssymb}
\usepackage{latexsym}
\usepackage{eufrak}

\theoremstyle{plain}
\newtheorem{theorem}[equation]{Theorem}
\newtheorem{proposition}[equation]{Proposition}
\newtheorem{lemma}[equation]{Lemma}

\theoremstyle{definition}

\theoremstyle{remark}

\newtheorem{remark}[equation]{Remark}

\makeatletter
\renewcommand{\subsection}{\@startsection{subsection}{2}{0pt}{-3ex
plus -1ex minus -0.2ex}{-2mm plus -0pt minus
-2pt}{\normalfont\bfseries}} \makeatother

\numberwithin{equation}{subsection}

\DeclareMathOperator{\End}{\mathrm{End}}

\DeclareMathOperator{\reg}{{\mathrm{reg}}}

\DeclareMathOperator{\Tr}{\mathrm{Tr}}

\newcommand{\dis}{\displaystyle}
\newcommand{\erem}{\hfill$\lozenge$\end{remark}}
\newcommand{\beq}{\begin{equation}\label}
\newcommand{\eeq}{\end{equation}}

\newcommand{\f}[1]{\mathfrak{#1}}

\DeclareMathOperator{\Spec}{\mathrm{Spec}}

\newcommand{\iso}{{\,\stackrel{_\sim}{\longrightarrow}\,}}

\DeclareMathOperator{\Proj}{\mathrm{Proj}}

\renewcommand{\mid}{\enspace\big|\enspace}

\newcommand{\too}{\,\,\longrightarrow\,\,}

\newcommand{\mto}{\longmapsto}
\newcommand{\map}{\longrightarrow}
\newcommand{\onto}{\,\twoheadrightarrow\,}
\newcommand{\cd}{\!\cdot\!}

\newcommand{\vi}{${\sf {(i)}}\;$}
\newcommand{\vii}{${\sf {(ii)}}\;$}

\newcommand{\g}{{\mathfrak{g}}}

\newcommand{\M}{{\mathscr{M}}}

\newcommand{\inv}{^{-1}}
\newcommand{\oo}{{\mathcal O}}

\renewcommand{\part}{{\f P}}

\newcommand{\ff}{{\mathbf{f}}}
\newcommand{\bv}{{\mathbf{v}}}
\newcommand{\x}{{\mathbf{x}}}
\newcommand{\y}{{\mathbf{y}}}
\renewcommand{\c}{\C[\x,\y]}

\def\C{{\mathbb{C}}}

\def\gln{{\mathfrak{g}\mathfrak{l}}_n(\C)}

\newcommand{\eps}{\epsilon}

\newcommand{\h}{{\mathfrak{h}}}
\newcommand{\sset}{\subset}

\newcommand{\en}{\enspace}

\newcommand{\into}{{\,\hookrightarrow\,}}

\newcommand{\Z}{{\mathbb{Z}}}

\newcommand{\Hilb}{{\operatorname{Hilb}^n{\mathbb{\C}}^2}}

\begin{document}
\centerline{\Large{\textbf{\Large{A remark on a theorem of M. Haiman}}}}

\bigskip
\centerline{\sc 
Victor Ginzburg}
\medskip
\begin{abstract} We deduce
a special case of a theorem of M. Haiman concerning alternating
polynomials in $2n$ variables from
 our results about {\em almost commuting variety},
obtained  earlier in a joint work with W.-L. Gan.
\end{abstract}
\section{Introduction}
\subsection{Main result.}\label{main} Write 
 $\C[\x,\y]:=
\C[x_1,\ldots,x_n,y_1,\ldots,y_n]$
for a polynomial ring in two sets of variables
 $\x=(x_1,\ldots,x_n)$ and $\y=(y_1,\ldots,y_n)$.
The  Symmetric group $S_n$ acts naturally on the $n$-tuples
$\x$ and $\y$,
and this gives rise to an $S_n$-diagonal action
on the algebra $\C[\x,\y].$
We write $\c^{S_n}\sset \c$ for the subalgebra
of $S_n$-invariant polynomials and $A:=\c^\eps\sset\c$
for the subspace of $S_n$-alternating
polynomials. 
The space $A$
is stable under multiplication by elements of the algebra
$\c^{S_n}$, in particular, it may be viewed as a module
over $\C[\y]^{S_n}\sset\c^{S_n}$, 
 the subalgebra
of symmetric polynomials in the last $n$ variables $y_1,\ldots,y_n$.

For each $k=1,2,\ldots,$ let
$A^k$ be the $\C$-linear subspace in $\c$
spanned by the products of $k$ elements of $A$. The action of
$\C[\y]^{S_n}$
on $A$ induces one on $A^k$, hence each  space $A^k,\,k=1,2,\ldots,$
acquires a natural $\C[\y]^{S_n}$-module structure.

The goal of this note is to give a direct proof of 
 the following  special case of
a much stronger result due to M. Haiman  \cite[Proposition 3.8.1]{Ha2}. 

\begin{theorem}\label{haiman} For each  $k=1,2,\ldots,$ 
the space $A^k$ is a free 
 $\C[\y]^{S_n}$-module.\end{theorem}

In  an earlier paper, Haiman showed, cf.
\cite[Proposition 2.13]{Ha1}, that  the above
theorem holds for all $k\gg0$.
The corresponding statement for {\em all}
$k$ follows from Haiman's proof 
of his
 his Polygraph theorem, the main technical result in \cite{Ha2}.

\subsection{Reminder from \cite{GG}.}\label{int1}
Let $V:=\C^n$ and let  $\g := \End(V)=\gln$ 
be the Lie algebra of $n\times n$-matrices.
We will write 
elements of $V$ as column vectors, and elements of
$V^*$ as row vectors.
Following Nakajima \cite{Na},
we consider the following
closed affine subscheme in the vector space
$\g\times\g\times V\times V^*$:
\begin{equation}\label{M}
 \M := \{ (X,Y,i,j)\in \g\times\g\times
V\times V^* \mid [X,Y]+ij=0\}. 
\end{equation}
More precisely, let 
$\C[\g\times\g\times V\times V^*]=\C[X,Y,i,j]$ denote the polynomial
algebra, and let $J\subset \C[X,Y,i,j]$ be the ideal
generated by the $n^2$ entries of the matrix $[X,Y]+ij$.
Then, by definition, we have $\M=\Spec\C[X,Y,i,j]/J,$
a not necessarily reduced affine scheme.

To describe the structure of the scheme $\M$,
for each integer $k\in\{0,1,\ldots,n\}$, set
$$ \M'_k := \bigg\{ (X,Y,i,j)\in\M \left |
\begin{array}{c}
\textrm{$Y$ has pairwise distinct eigenvalues},\\
\dim(\C[X,Y]i)=n-k,\quad \dim(j\C[X,Y])=k \end{array}\right\} $$
and let $\M_k$ be the closure of $\M'_k$ in $\M$.

The following result was proved in 
\cite[Theorem 1.1.2]{GG}.

\begin{theorem} \label{t1}
\vi The irreducible components of $\M$ are $\M_0, \ldots, \M_n$.

\vii $\M$ is a reduced  complete intersection in $\g\times\g\times
V\times V^*$.
\end{theorem}

\subsection{Geometric interpretation of $A^k$.} 
Write
$\h:=\C^n$ for the
tautological permutation representation of  the Symmetric group $S_n$,
and let $S_n$ act diagonally on
$\h\times\h$. 
We have an obvious identification
$\C[\h\times\h]=\c$, in particular, we may view 
the vector space $A^k$, see \S\ref{main}, as a subspace in $\C[\h\times\h].$

The quotient $(\h\times\h)/S_n$ has a natural
structure of algebraic variety, with
coordinate ring 
$$\C[(\h\times\h)/S_n]=\C[\h\times\h]^{S_n}=\c^{S_n}.$$

We may also view $\h$ as the Cartan subalgebra in $\g$ formed by diagonal
matrices, so
we have a tautological imbedding $\h\times\h\into\g\times\g$.
We define the following map
\beq{j}\jmath:\
 \h\times\h\into \g\times\g\times V\times V^*,
\quad(\x,\y)\mto(\x,\y,i_o,0),
\eeq
where $i_o$ stands for the vector $i_o:=(1,1,\ldots,1)\in V$. 

The group $G=GL(V)$ acts naturally on $\M$ by the formula
$g: (X,Y,i,j)\mto$
$\dis (gXg\inv,gYg\inv, g\cdot i, j\cdot g\inv).$
Note that
$S_n$, viewed as the subgroup of permutation matrices in $G$,
fixes the vector $i_o\in V$. So, the image of $\jmath$ is an $S_n$-stable
subset in $\g\times\g\times V\times V^*$.
Thus, the map $\jmath$ gives
an $S_n$-equivariant closed imbedding $\jmath: \h\times\h\into \M$.

The action of $G$ on $\M$ gives rise to a $G$-action $g: f\mapsto g(f)$ 
on the coordinate ring $\C[\M]$ by algebra
automorphisms.
For each $k=1,2,\ldots,$ we set
$$ \C[\M]^{(k)}:=\{f\in\C[\M]\mid g(f)=(\det g)^k \cd f,\quad\forall
g\in G\}.
$$

A key ingredient in our approach to Theorem \ref{haiman} is
the following result to be proved in the next section.
\begin{proposition}\label{A^k} For each $k=1,2,\ldots,$
restriction of functions  via $\jmath$
induces a vector space isomorphism $\jmath^*: \C[\M]^{(k)}\iso A^k$.
\end{proposition}

\begin{remark} According to \cite[Lemma 2.9.2]{GG},
 restriction of functions via $\jmath$
induces also a  graded algebra
isomorphism
$$
\jmath^*:\
\C[\M]^G\iso
\C[\h\times\h]^{S_n}.
$$
The latter isomorphism may be viewed as a version
of Proposition \ref{A^k} for $k=0.$
\end{remark}

\begin{remark}
Let  $\Hilb$ denote the Hilbert scheme
of $n$ points in the plane, cf. e.g. \cite{Ha1},\cite{Na}.
The Hilbert scheme comes equipped with a natural
ample line bundle $\oo(1)$, cf. \cite{Ha1}.

We remind the reader that, for each $k=1,2,\ldots,$
M. Haiman constructed in \cite{Ha1} a natural map $A^k\to
\Gamma(\Hilb,\oo(k))$;
Moreover, it follows from the results of \cite{Ha2}
that this map is an isomorphism.
\end{remark}
\section{Proof of Proposition \ref{A^k}}
\subsection{} Fix nonzero volume elements
$\bv\in \wedge^n V$ and $\bv^*\in \wedge^n V^*,$ respectively.
Given an $n$-tuple $\ff=(f_1,\ldots,f_n)\in\C\langle x,y\rangle,$
of  noncommutative polynomials in two variables, we
consider polynomial functions $\psi,\phi\in\C[\g\times\g\times
V\times V^*]$
of the form 
\begin{align}\label{bv}
&\psi_\ff(X,Y,i,j)=\langle\bv^*,\, f_1(X,Y)i\wedge\ldots
f_n(X,Y)i\rangle,\\
 &\phi_\ff(X,Y,i,j)=\langle
jf_1(X,Y)
\wedge\ldots
jf_n(X,Y),\,\bv\rangle,\nonumber
\end{align}
where $f_r(X,Y)$ denotes the matrix obtained by plugging the two
matrices $X,Y\in\g$ in the  noncommutative polynomial $f(x,y)$.
We will keep the notation
$\psi_\ff,\phi_\ff$ for the restriction of the corresponding
function to the closed subvariety $\M\sset \g\times\g\times
V\times V^*.$ It is clear that, restricting these functions further
to the subset $\h\times\h\sset\M$, on has
$\jmath^*\psi_\ff\in A$ and $\jmath^*\phi_\ff=0$. 

Recall that $\M=\M_0\cup\ldots\cup \M_n$, is a union of $n+1$
irreducible components. It is immediate from the
definition
of the set $\M_r$, cf. \S\ref{int1}, that, for any choice of 
 $n$-tuple $\ff=(f_1,\ldots,f_n),$ the function
$\psi_\ff$ vanishes on $\M_r$ whenever $r\neq 0$, while
$\phi_\ff$ vanishes on $\M_r$ whenever $r\neq n$.
Since each irredicible component is {\em reduced}, by Theorem \ref{t1},
the above vanishings hold
scheme-theoretically:
\beq{van}
\psi_\ff|_{\M_r}=0\quad\forall r\neq0,
\quad\text{and}\quad\phi_\ff|_{\M_r}=0\quad\forall r\neq n.
\eeq

Next, similarly to $\C[\M]^{(k)},$ for each $k\in\Z$, we introduce the space
 $\C[\g\times\g\times V\times V^*]^{(k)}$ 
of polynomial functions on $\g\times\g\times V\times V^*$
that satisfy the equation $g(f)=(\det g)^k \cd f,\quad\forall
g\in G.$ It is clear that 
$$\psi_\ff\in \C[\g\times\g\times V\times V^*]^{(1)},
\quad\text{resp.},\quad
\phi_\ff\in \C[\g\times\g\times V\times V^*]^{(-1)},
\quad\forall \ff=(f_1,\ldots,f_n).
$$
Observe  that
$\C[\g\times\g\times V\times V^*]^{(k)}$ 
is naturally a $\C[\g\times\g\times V\times V^*]^G$-module.
Applying Weyl's fundamental theorem on $GL_n$-invariants
we deduce that this $\C[\g\times\g\times V\times V^*]^G$-module
is generated by  products  of the form
$\psi_1\cdot\ldots\cdot\psi_p\cdot\phi_1\cdot\ldots\cdot\phi_q$,
where $p-q=k$ and where each factor $\psi_r$, resp. each factor
 $\phi_r$, is of the form $\psi_\ff$, resp., $\phi_\ff$.

The action of $G$ on $\C[\g\times\g\times V\times V^*]$
being completely reducible, we deduce that restricting
functions from $\g\times\g\times V\times V^*$ to $\M$
yields a surjection
$\C[\g\times\g\times V\times V^*]^{(k)}\onto
\C[\M]^{(k)},$ for any $k\in\Z$. It follows that
$\C[\M]^{(k)},$ viewed as a $\C[\M]^G$-module,
is again generated by the products
$\psi_1\cdot\ldots\cdot\psi_p\cdot\phi_1\cdot\ldots\cdot\phi_q$,
with $p-q=k$. Furthermore, from \eqref{van}, we see that 
for  $k\geq 0$
we must have $p=k\en\&\en q=0$. On the other hand, for  $k\leq 0$
we must have $p=0\en\&\en q=k$.

From now on, we assume that $k\geq 1$. Thus,
 the imbedding $\M_0\into \M$ induces a bijection
$\C[\M]^{(k)}\iso \C[\M_0]^{(k)}.$ It follows that
$\C[\M]^{(k)}$  is generated, as a $\C[\M]^G$-module, by the products
$\psi_1\cdot\ldots\cdot\psi_k$.
 Since 
 $\jmath^*\psi_\ff\in A$ for any $\ff$, we find 
 that $\jmath^*(\psi_1\cdot\ldots\cdot\psi_k)\in A^k$, hence
$\jmath^*(\C[\M]^{(k)})=\jmath^*(\C[\M_0]^{(k)})\sset A^k.$

To prove injectivity of the restriction map
$\jmath^*: \C[\M]^{(k)}\iso A^k,$ we observe that $G\cdot
\jmath(\h\times\h),$ the $G$-saturation of the image of the imbedding $\jmath$,
is an irreducible variety of dimension $n^2+n=\dim \M$. Furthermore, for
any diagonal matrix $Y\in\h$ with pairwise distinct eigenvalues,
we have $\C[Y]i_o=V$. Hence, a Zariski open subset of 
$G\cdot\jmath(\h\times\h)$ is contained $\M'_0$, cf. \S\ref{int1}. 
Since $\M_0=\overline{\M'_0}$ and $G\cdot\jmath(\h\times\h)$ is 
irreducible, we conclude that $G\cdot\jmath(\h\times\h)\sset \M_0$
and, moreover, the set $G\cdot\jmath(\h\times\h)$ is Zariski dense
in $\M_0$. Thus, for any $f\in \C[\M_0]^{(k)}$ such that $\jmath^*(f)=0$
we must have $f=0$. This proves injectivity of the map $\jmath^*$.

Observe next that since $\C[\M]^{(k)}$  is generated, as a $\C[\M]^G$-module, by the products
$\psi_1\cdot\ldots\cdot\psi_k$ it suffices to prove surjectivity
 of the map $\jmath^*$ for $k=1$. To prove the latter, we identify
$A=\c^\epsilon$ with $\wedge^n\C[x,y]$, the $n$-th exterior power
of the vector space $\C[x,y]$ of polynomials in 2 variables.
With this identification, the space $A$ is spanned by expressions
of the form $f_1\wedge\ldots\wedge f_n,\,f_1,\ldots,f_n\in\C[x,y]$.

Now, recall that for any $(X,Y,i,j)\in\M_0$ we have
$[X,Y]=[X,Y]+ij=0$. Therefore, for any $f\in\C[x,y]$,
the expression $f(X,Y)$ is a well-defined matrix. In other words, for
any lift of $f$ to a noncommutative polynomial $\hat{f}\in\C\langle
x,y\rangle,$
i.e., for any $\hat{f}$ in the preimage of $f$ under the natural projection
$\C\langle
x,y\rangle\onto\C[x,y]$, we have $\hat{f}(X,Y)=f(X,Y).$
Thus, given an $n$-tuple $f_1,\ldots,f_n\in\C[x,y]$,
we have a well-defined element
$$\psi_\ff=\langle\bv^*, f_1(X,Y)i\wedge\ldots
f_n(X,Y)i\rangle\in \C[\M]^{(1)}.
$$
It is straightforward to verify that
$\jmath^*\psi_\ff=f_1\wedge\ldots\wedge f_n.$
This proves 
 surjectivity of the map $\jmath^*$ 
and completes the proof of the Proposition.
\newpage

\section{Proof of Theorem \ref{haiman}.} 
\subsection{A flat morphism.} \label{flat_s}
Denote by $\C^{(n)}$ the set of unordered $n$-tuples of complex
numbers. Let 
$$
\pi:\ \M\too \C^{(n)},\quad (X,Y,i,j)\mto\Spec Y
$$
 be the map that sends
$(X,Y,i,j)$ to the unordered $n$-tuple $\Spec Y$ of eigenvalues of $Y$,
counted with multiplicities. 

According to \cite[Proposition 2.8.2]{GG}, the morphism $\pi$ is {\em
flat}. In algebraic terms, this means
that $\C[\M]$ is a flat $\C[\y]^{S_n}$-module.

\subsection{}
We have the  standard grading
$\dis\C[\y]^{S_n}=\oplus_{d\geq 0}\C^d[\y]^{S_n}$,
 by degree of the polynomial.
Write $\dis\C[\y]^{S_n}_+:=\oplus_{d> 0}\C^d[\y]^{S_n}$ 
for the augmentation ideal.

Let $E$ be any  flat
nonnegatively graded  $\C[\y]^{S_n}$-module.
Then, choosing representatives in $E$ of a
 $\C$-basis of the vector space
$\dis E/\C[\y]^{S_n}_+E$ yields a {\em free} $\C[\y]^{S_n}$-basis in $E$.
Hence, any  flat
nonnegatively graded  $\C[\y]^{S_n}$-module is {\em free}.

Next, we have a  $\C^\times$-action on $\M$ given
by the formula
 $\dis\C^\times\ni z: (X,Y,i,j)\mto (z\cdot X,z\cdot Y,z\cdot i,z\cdot j)$.
This  $\C^\times$-action gives rise to a natural grading 
$\dis\C[\M]=\oplus_{d\geq 0}\C^d[\M]$, on the algebra
$\C[\M]$.
With this grading, the pull-back morphism $\pi^*: \C[\y]^{S_n}
\map \C[\M]$ induced by the map $\pi: (X,Y,i,j)\mto \Spec Y,$
is a graded algebra morphism, so $\C[\M]$ may be viewed as
a {\em graded} $\C[\y]^{S_n}$-module.

Thus, according to \S\ref{flat_s} we know
that $\C[\M]$ is a flat nonnegatively graded  $\C[\y]^{S_n}$-module.
As has been explained at the beginning of the proof,
this implies that $\C[\M]$ is  free over $\C[\y]^{S_n}$.
Further, by complete reducibility of the $G$-action on $\C[\M]$,
we deduce that $\C[\M]^{(k)}$, viewed as a graded $\C[\y]^{S_n}$-submodule
in $\C[\M]$, splits off as a direct summand,
for any $k=1,2,\ldots$. Hence,
$\C[\M]^{(k)}$ is projective, in particular, flat over $\C[\y]^{S_n}$,
 as a direct summand of  a free $\C[\y]^{S_n}$-module.

To complete the proof of 
Theorem \ref{haiman} we observe
 that $\C[\M]^{(k)}$ is a graded
 $\C[\y]^{S_n}$-module. Thus,  we conclude as above that 
this graded module must be free
 over $\C[\y]^{S_n}$.

\setcounter{equation}{0}
\footnotesize{

\noindent                                                                          
{\bf V.G.}: 
Department of Mathematics, University of Chicago,
Chicago, IL 60637, USA;\\
\hphantom{x}\quad\, {\tt ginzburg@math.uchicago.edu}}


\begin{thebibliography}{APK}










\bibitem[GG]{GG} W.-L. Gan, V. Ginzburg,
 {\em Almost-commuting variety, D-modules, and Cherednik Algebras.}\hfill\break
  {\tt{arXiv:math.RT/0409262}}
{\em }
\bibitem[Ha1]{Ha1} M. Haiman, $t,q$-{\em Catalan numbers and the Hilbert 
scheme.} Selected papers in honor of Adriano Garsia (Taormina, 1994).  
Discrete Math.  \textbf{193}  (1998),  no. 1-3, 201--224. 

\bibitem[Ha2]{Ha2} M. Haiman,
{\em
Hilbert schemes, polygraphs and the Macdonald positivity conjecture.}
  J. Amer. Math. Soc. \textbf{14}  (2001),   941--1006.


\bibitem[Na]{Na} H. Nakajima,
{\em  Lectures on Hilbert schemes of points on surfaces},
University Lecture Series, 18,
Amer. Math. Soc., Providence, RI, 1999.

\end{thebibliography}
\end{document}